\documentstyle[11pt]{article}
\onecolumn

\begin{document}

\title{Symmetric Logic Synthesis with Phase Assignment}

\author{ {\sc Nico F. Benschop} }

\date{ {\it 22nd Information Theory Symposium (IEEE}/Benelux, U-Twente, May 2001)}

\maketitle

\begin{abstract} \small{
Decomposition of any Boolean Function $BF_n$ of $n$ binary inputs
into an optimal inverter coupled network of Symmetric Boolean
functions $SF_k~(k \leq n)$ ~is described. Each $SF$ component
is implemented by Threshold Logic Cells, forming a complete and
compact $T$-$Cell$ Library. Optimal phase assignment of input
polarities maximizes local symmetries. $Rank~spectrum$
is a new $BF_n$ description independent of input ordering,
obtained by mapping its minterms onto an othogonal $n \times n$
grid of (transistor-) switched conductive paths, minimizing
crossings in the silicon plane. Using this ortho-grid structure
for the layout of $SF_k$ cells, without mapping to $T$-cells,
yields better area efficiency, exploiting the maximal logic
path sharing in $SF$'s.
Results obtained with a $CAD$ tool "$Ortolog$" based on these
concepts, are reported. Relaxing symmetric- to $planar$- Boolean
functions is sketched, to improve low- symmetry $BF$ decomposition.
}
\end{abstract}

\section{Introduction}    

Since the early eighties the synthesis of combinational logic for the
design of integrated circuits ($IC's$) is increasingly automated.
Present logic synthesis tools, near the bottom of the $IC$ design
hierarchy, just above layout, is fairly mature, being intensively
applied in the design of production IC's. But some problems remain:

{\bf A}.~
  Logic synthesis tools often have a disturbing {\bf order dependence}.
Re-ordering signals, which should not affect the result, can cause a
considerable increase or decrease of silicon area. To curb computer
time, synthesis tools avoid global analysis which tends to grow
exponentially with the number of inputs. Hence a local approach is
preferred, using a greedy algorithm, taking the first improvement
that comes along. The result then depends on the ordering of cubes in
a $PLA$ listing, or the input order in a $BDD$ (binary decision diagram)
[1][2][3] representing a Boolean function $(BF)$. This effect is reduced
by {\bf global} analysis, and by symmetric function components $SF$,
being independent of input ordering. $CPU$ time is reduced by the
'arithmetization' via $spectral~ BF_n$ analysis, a new method of
characterizing $BF$'s, to be explained.

{\bf B}.~ Optimal polarity or {\bf phase assigment} of signals,
either inputs or intermediate variables, is still an unsolved problem,
although some heuristics are applied. Input phases influence logic
symmetries, to be exploited for an efficient decomposition, that is
essentially synthesis.

{\bf C}.~
The use of a standard {\bf cell library} is forcing decomposition-
and cell mapping stages to produce a sub optimal gate network,
versus {\bf compiled cells} as needed [4]: using no cell library
but a programmable grid template, to be discussed.
The proposed 'orthogrid' $BF$ structure is an experiment in that
direction, to be extended to planar $F$'s beyond symmetric $F$'s
as 'grid template' alternative to $FPGA$ or $FPMUX$ cells [5].
Performance prediction, that comes with a cell library, is then
done by the {\bf cell compiler}, which is quite feasible,
replacing library maintenance by compiler support.

{\bf D}.~
{\bf Complete testing} of combinational logic circuits requires
{\bf irredundancy}, guaranteed only in $sum~of~cubes$ 2-level
implementation. Logic in factored form, the usual result of a
synthesis tool, sometimes has testability problems. Restriction
to a {\bf disjoint product} is proposed, with factors having no
common input. This guarantees the irredundancy needed for $BF$
testability in factored form. And: disjoint products yield a
$spectral ~calculus$,with a $BF ~rank~spectrum$ independent of
input ordering, and a convolution composition rule.

{\bf Order independent Logic Synthesis}:~      

The mentioned problems in present synthesis $CAD$ imply that no
optimality (nor full testability) is guaranteed, nor does one know
how close/far the optimum is. Presently, only by many synthesis runs
(design space exploration) a feeling is obtained for the complexity
of the functions to be synthesized, allowing a trade-off between
circuit area , -delay, and power dissipation, however at a high
cost in $CPU$ time.

Our $aim$ is to improve this situation, crucial for the future of
digital $VLSI$ systems. The emphasis is on order- independent function
representation, using a spectral technique called $rank~spectrum$,
and on $global$ analysis before synthesis, which then becomes
feasible. In fact we go one step beyond $BDD$ type of $BF$
descriptions, by mapping minterms as paths in an orthogonal grid,
using  symmetric $F$'s and signal phasing.

Then methods similar to those applied in signal processing, like
the frequency spectrum, or convolution of impulse response and
input sequence in the time domain, can also be applied to Boolean
functions. This yields:\\
--- synthesis by global structure analysis,\\
--- with arithmetization of Boolean algebra \\
--- via a rank-spectrum technique.

\section{Ortho grid, rank spectrum}    

{\it Def}:
   {\bf Orthogrid plot}: map each $minterm$ of a $BF_n$ (as 0/1
string of length $n$) in an orthogonal grid, as an $n$-step path
from the origin to the $n$-th diagonal. In input sequence,
step down if '0', and right if '1' (see fig. 1).\\This models a
pass transistor network on silicon, with a conducting path from
the origin to the $n$-th diagonal for the given minterm. $OR$-ing
all paths yields function $F$=1 only if some path connects
origin to final diagonal.

For $n$ inputs, each path ends on the $n$-th diagonal. All minterms
of equal {\bf rank} (number of ones) end in the same point on the
$n$-th diagonal. Without confusion such minterm-set is also called
a $rank$ of $F$. For the orthogrid plot of a single rank
$XOR$ product (4 terms, rank 2) see fig.1.

$Def$:~ a {\bf rank function} $RF$  has only one non-empty
     rank (equal rank minterms).

$Def$:~ $BF_n$ {\bf rank spectrum} is the vector
     of path (minterm) counts per rank [0 - $n$]

A $BF_n$ is the sum of its rank functions, and its rank spectrum
is independent of input ordering. In general, crossing paths are not
allowed to touch each other, to be drawn with a bridge or tunnel.
This makes the ~$ortho~grid$ style cumbersome for larger functions,
and probably explains the popularity of the Shannon-tree, which can
be displayed free of crossings, that is: as a planar a-cyclic graph.
However, {\bf path sharing} is essential to recognize {\bf common
factors}, which is a clue to logic synthesis, showing the power of
$BDD$'s and the othogrid representation.

{\bf Planar node: factoring paths}       

$Def$: a node is {\bf planar} if all paths connect there.
~~(e.g. the circled node in fig.1).\\ So all such paths are cut in two
parts: each first section from the origin is continued (multiplied)
by all second sections to the final diagonal.\\ A function $F$ with
all paths (minterms) passing through a planar node is a product of
two functions $F=G*H$ sharing no inputs, where $G$ is a rank
function; here $G(a,b)$ and $H(c,d)$. A planar node plays the role
of a $factor~node$. {\it Planarization} is essential for synthesis,
obtained by proper choice of order and polarity of inputs.

%
\setlength{\unitlength}{0.01in}%
\begingroup\makeatletter\ifx\SetFigFont\undefined
\def\x#1#2#3#4#5#6#7\relax{\def\x{#1#2#3#4#5#6}}%
\expandafter\x\fmtname xxxxxx\relax \def\y{splain}%
\ifx\x\y   
\gdef\SetFigFont#1#2#3{%
  \ifnum #1<17\tiny\else \ifnum #1<20\small\else
  \ifnum #1<24\normalsize\else \ifnum #1<29\large\else
  \ifnum #1<34\Large\else \ifnum #1<41\LARGE\else
     \huge\fi\fi\fi\fi\fi\fi
  \csname #3\endcsname}%
\else
\gdef\SetFigFont#1#2#3{\begingroup
  \count@#1\relax \ifnum 25<\count@\count@25\fi
  \def\x{\endgroup\@setsize\SetFigFont{#2pt}}%
  \expandafter\x
    \csname \romannumeral\the\count@ pt\expandafter\endcsname
    \csname @\romannumeral\the\count@ pt\endcsname
  \csname #3\endcsname}%
\fi
\fi\endgroup
\begin{picture}(700,210)(80,90)
\put(150,290){\vector(1,0){85}}
\put(200,145){\vector(1,1){85}}
\put( 85,230){\vector(0,-1){85}}
\put(120,220){\framebox(40,40){}}
\put(160,180){\framebox(40,40){}}
\multiput(120,180)( 7.27273, 0.00000){6}{\line( 1, 0){3.636}}
\multiput(200,220)( 0.00000, 7.27273){6}{\line( 0, 1){3.636}}
\multiput(205,140)(15.51111,20.68148){4}{\line( 1, 1){8.667}}
\multiput(120,260)( 0.00000,-8.00000){20}{\line( 0,-1){4.000}}
\multiput(275,260)(-8.00000, 0.00000){20}{\line(-1, 0){4.000}}
\multiput(120, 95)(21.39803,21.39803){ 8}{\line( 1, 1){10.214}}
\put(300,200){\shortstack
     {{\it a b c d}\\---------\\1 0 1 0\\1 0 0 1\\0 1 1 0\\0 1 0 1}}
\put(300,160){\shortstack{N = 6\\L = 8}}
\put( 80,120){\makebox(0,0)[lb]{\smash{\SetFigFont{14}{14.4}{bf}0}}}
\put(240,285){\makebox(0,0)[lb]{\smash{\SetFigFont{14}{14.4}{bf}1}}}
\put( 90,275){\makebox(0,0)[lb]{\smash{\SetFigFont{14}{14.4}{it}F}}}
\put(140,270){\makebox(0,0)[lb]{\smash{\SetFigFont{14}{14.4}{it}a}}}
\put(180,270){\makebox(0,0)[lb]{\smash{\SetFigFont{14}{14.4}{it}b}}}
\put(220,270){\makebox(0,0)[lb]{\smash{\SetFigFont{14}{14.4}{it}c}}}
\put(260,270){\makebox(0,0)[lb]{\smash{\SetFigFont{14}{14.4}{it}d}}}
\put(100,245){\makebox(0,0)[lb]{\smash{\SetFigFont{14}{14.4}{it}\_}}}
\put(100,230){\makebox(0,0)[lb]{\smash{\SetFigFont{14}{14.4}{it}a}}}
\put(100,200){\makebox(0,0)[lb]{\smash{\SetFigFont{14}{14.4}{it}\_}}}
\put(100,185){\makebox(0,0)[lb]{\smash{\SetFigFont{14}{14.4}{it}b}}}
\put(100,155){\makebox(0,0)[lb]{\smash{\SetFigFont{14}{14.4}{it}\_}}}
\put(100,140){\makebox(0,0)[lb]{\smash{\SetFigFont{14}{14.4}{it}c}}}
\put(100,110){\makebox(0,0)[lb]{\smash{\SetFigFont{14}{14.4}{it}\_}}}
\put(100, 95){\makebox(0,0)[lb]{\smash{\SetFigFont{14}{14.4}{it}d}}}
\thinlines
\put(160,220){\circle{10}}
\put(135, 95){\makebox(0,0)[lb]{\smash{\SetFigFont{14}{14.4}{rm}0}}}
\put(170,135){\makebox(0,0)[lb]{\smash{\SetFigFont{14}{14.4}{rm}1}}}
\put(205,175){\makebox(0,0)[lb]{\smash{\SetFigFont{14}{14.4}{rm}2}}}
\put(240,210){\makebox(0,0)[lb]{\smash{\SetFigFont{14}{14.4}{rm}3}}}
\put(275,245){\makebox(0,0)[lb]{\smash{\SetFigFont{14}{14.4}{rm}4}}}
\put(220,155){\makebox(0,0)[lb]{\smash{\SetFigFont{14}{14.4}{it}rank}}}
\put(170,110){\makebox(0,0)[lb]{\smash{\SetFigFont{14}{14.4}{it}spectrum =}}}
\put(270,110){\makebox(0,0)[lb]{\smash{\SetFigFont{14}{14.4}{rm}[0,0,4,0,0]}}}
\put(380,110){\makebox(0,0)[lb]{\smash{\SetFigFont{14}{14.4}{it}F =(a\#b)(c\#d)}}}
%
\thicklines
\put(420,210){\framebox(40,40){}}
\put(420,230){\line(1,0){ 40}}
\put(440,235){\line(0,1){ 15}}
\put(440,235){\line(1,0){  5}}
\put(445,225){\line(0,1){ 10}}
\put(440,225){\line(1,0){  5}}
\put(440,210){\line(0,1){ 15}}
\put(475,210){\shortstack{N = 9\\L = 12}}
\thinlines
\multiput(460,250)(7.27273, 0.00000){6}{\line( 1, 0){  3.636}}
\multiput(420,210)(0.00000,-7.27273){6}{\line( 0,-1){  3.636}}
\put(400,150)
  {\makebox(0,0)[lb]{\smash{\SetFigFont{12}{12.0}{rm}spectr.[0 0 4 0 0]}}}
\put(425,255){\makebox(0,0)[lb]{\smash{\SetFigFont{14}{14.4}{it}a}}}
\put(450,255){\makebox(0,0)[lb]{\smash{\SetFigFont{14}{14.4}{it}c}}}
\put(475,255){\makebox(0,0)[lb]{\smash{\SetFigFont{14}{14.4}{it}b}}}
\put(495,255){\makebox(0,0)[lb]{\smash{\SetFigFont{14}{14.4}{it}d}}}
\put(580,250){\line(0,-1){ 80}}
\put(580,250){\line(1, 0){ 85}}
\put(580,210){\line(1, 0){ 40}}
\put(620,210){\line(0, 1){ 40}}
\put(635,210){\shortstack{N = 11\\L = 12}}
\put(560,150)
  {\makebox(0,0)[lb]{\smash{\SetFigFont{12}{12.0}{rm}spectr.[1 0 2 0 1]}}}
\put(590,255){\makebox(0,0)[lb]{\smash{\SetFigFont{14}{14.4}{it}a}}}
\put(595,247){\makebox(0,0)[lb]{\smash{\SetFigFont{14}{14.4}{it}.}}}
\put(610,265){\makebox(0,0)[lb]{\smash{\SetFigFont{14}{14.4}{it}\_}}}
\put(610,255){\makebox(0,0)[lb]{\smash{\SetFigFont{14}{14.4}{it}b}}}
\put(635,265){\makebox(0,0)[lb]{\smash{\SetFigFont{14}{14.4}{it}\_}}}
\put(635,255){\makebox(0,0)[lb]{\smash{\SetFigFont{14}{14.4}{it}c}}}
\put(640,247){\makebox(0,0)[lb]{\smash{\SetFigFont{14}{14.4}{it}.}}}
\put(655,255){\makebox(0,0)[lb]{\smash{\SetFigFont{14}{14.4}{it}d}}}
\put(617,225){\makebox(0,0)[lb]{\smash{\SetFigFont{14}{14.4}{rm}-}}}
\put(577,225){\makebox(0,0)[lb]{\smash{\SetFigFont{14}{14.4}{rm}-}}}
\put(577,185){\makebox(0,0)[lb]{\smash{\SetFigFont{14}{14.4}{rm}-}}}
\put(595,207){\makebox(0,0)[lb]{\smash{\SetFigFont{14}{14.4}{rm}.}}}
\end{picture}

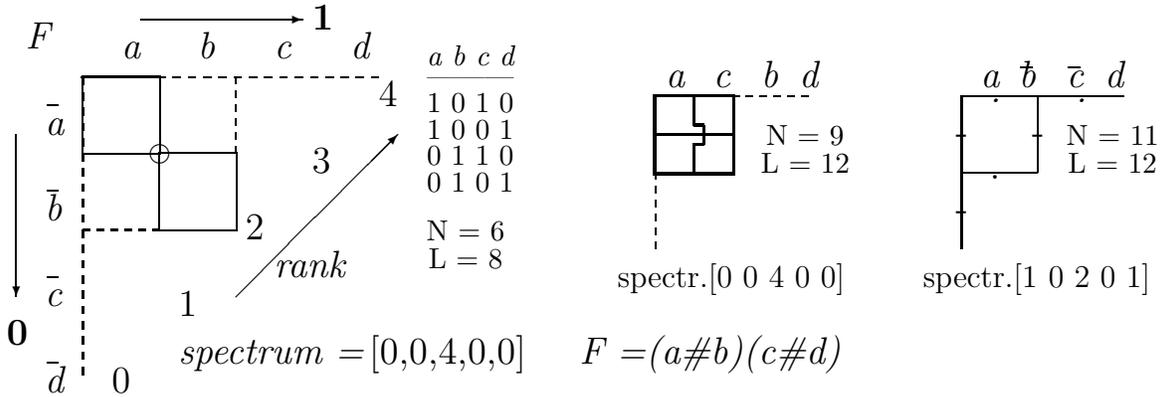
{\bf Fig 1}. Gridplot of~ F= XOR pair product.~~
{\bf Fig 2}. Planarize: permute/invert inputs

Counting occupied gridpoints (nodes), multiple for non-planar nodes,
yields a good criterion for a logic optimization algorithm
($planarization$):

{\bf Factoring criterion}:
Permute and invert ($phase$) inputs to ~minimize node count $N$.

Alternatively, the number of links $L$, counting the transistors,
could be minimized. Node count $N$ dominates over link count for
practical technological reasons. A $bridge$ requires two via's to
another metal level, costing more than a transistor which is
simply a polysilicon line crossing (self-aligned) a diffusion path.
Permuting and inverting inputs, factored form fig.1 has minimal
$(N,L)=(6,8)$ of the three gridplots of $F$.

This orthogrid representation allows characterization of special types
of Boolean functions such as symmetric-, planar- and rank- functions,
to be considered next. Notice the maximally $2^n$ minterms are plotted
in a square grid of $n^2$ nodes, by virtue of dense {\it path sharing}
as partial factors.
Actually a half square suffices, up to diagonal $n$; the other half
plane could be used for the complement or dual of $F$ (as in $CMOS$).

\section{Symmetric and Threshold $BF's$}    

The well known Pascal Triangle, displayed in orthogonal grid fashion
(fig.3), gives in each node the number $R(i,j)$ of all paths connecting
that node to the origin. This is easily verified by its generation
rule: $R(i,j)=R(i-1,j) + R(i,j-1)$ is the sum of its predecessor
node path counts. Induction yields the path counting rank spectrum.

The $XOR$-product function $F$ (fig.1) is not symmetric in all
inputs, but it has two {\ partial} symmetries or input equivalences
(permute without changing $F$), written $a \cong b$ and $c \cong d$.
The $Ortolog$ algorithm (sect.5) detects and enhances such partial
symmetries.

\setlength{\unitlength}{0.01in}%
\begingroup\makeatletter\ifx\SetFigFont\undefined
\def\x#1#2#3#4#5#6#7\relax{\def\x{#1#2#3#4#5#6}}%
\expandafter\x\fmtname xxxxxx\relax \def\y{splain}%
\ifx\x\y   
\gdef\SetFigFont#1#2#3{%
  \ifnum #1<17\tiny\else \ifnum #1<20\small\else
  \ifnum #1<24\normalsize\else \ifnum #1<29\large\else
  \ifnum #1<34\Large\else \ifnum #1<41\LARGE\else
     \huge\fi\fi\fi\fi\fi\fi
  \csname #3\endcsname}%
\else
\gdef\SetFigFont#1#2#3{\begingroup
  \count@#1\relax \ifnum 25<\count@\count@25\fi
  \def\x{\endgroup\@setsize\SetFigFont{#2pt}}%
  \expandafter\x
    \csname \romannumeral\the\count@ pt\expandafter\endcsname
    \csname @\romannumeral\the\count@ pt\endcsname
  \csname #3\endcsname}%
\fi
\fi\endgroup
\begin{picture}(100,100)(20,130)
\put( 40,200){\makebox(0,0)[lb]{\smash{\SetFigFont{12}{12.0}
    {rm}1 - 1 - 1 - 1 - 1}}}
\put( 40,180){\makebox(0,0)[lb]{\smash{\SetFigFont{12}{12.0}
    {rm}1 - 2 - 3 - 4}}}
\put( 40,160){\makebox(0,0)[lb]{\smash{\SetFigFont{12}{12.0}
    {rm}1 - 3 - 6}}}
\put( 40,140){\makebox(0,0)[lb]{\smash{\SetFigFont{12}{12.0}{rm}1 - 4}}}
\put( 40,120){\makebox(0,0)[lb]{\smash{\SetFigFont{12}{12.0}{rm}1}}}
\thinlines
\put( 75,125){\vector(1,1){60}}
\put( 20,205){\makebox(0,0)[lb]{\smash{\SetFigFont{14}{14.4}{it}F}}}
\put(140,170)
   {\makebox(0,0)[lb]{\smash{\SetFigFont{12}{12.0}{it}Ranks 0 . . 4}}}
\put(100,130)
   {\makebox(0,0)[lb]{\smash{\SetFigFont{12}{12.0}{rm}F=1: spectr[1 4 6 4 1]}}}
\end{picture}

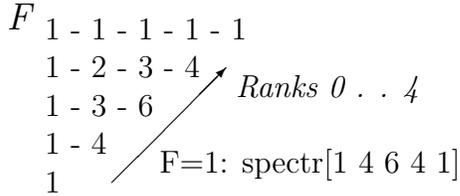
{\bf Fig 3}. ~Binomial path-count for full ranks.\\

A rank=2 symmetric function in 4 inputs contains all minterms of rank 2,
otherwise it cannot be an $SF$: there are $(4 ~choose ~2) = 6$ minterms,
in fact a {\bf full rank} has a {\bf binomial coefficient} number of
minterms. Notice in fig.1 there are two paths missing from a full
rank=2:~ 0011 and 1100 (see dotted lines), so $F$ is not symmetric.

{\bf --- Symmetric functions 'count' ---}       

{\it Def}: ~~a~ {\bf symmetric function} $SF$ \\ \hspace*{3mm}
does not change by permuting its inputs.

In other words, a function $SF$ is symmetric in all inputs if it
depends only on the number of 1-inputs, and not on their position.
Its ranks are either full or empty, so:

A symmetric function $SF[R]$ is determined
by the {\bf set} $R \subset [0, .. ,n]$ of its {\bf full ranks}.

An $n$-input function has $n$+1 ranks, with $2^{n+1}$ subsets, which
is the number of symmetric functions of $n$ inputs. For instance the
parity function is symmetric, active for an odd number of 1-inputs,
so the odd ranks are full, and all even ranks empty: $SF[odd]$.

\setlength{\unitlength}{0.01in}%
\begin{picture}(500,70)(0,0)
\thinlines
\put(320,20){\line(0,1){ 40}}
\put(320,20){\line(1,0){ 80}}
\put(340,20){\framebox(20,40){}}
\put(380,20){\framebox(20,40){}}
\put(450,20){\line(0,1){ 40}}
\put(450,20){\line(1,0){ 80}}
\put(490,20){\framebox(40,40){}}
\put(300,0){\makebox(0,0)[lb]{\smash{\SetFigFont{12}{12.0}{it}FA:}}}
\put(340,0){\makebox(0,0)[lb]{\smash{\SetFigFont{12}{12.0}{it}sum}}}
\put(380,0){\makebox(0,0)[lb]{\smash{\SetFigFont{12}{12.0}{rm}[1,3]}}}
\put(450,0){\makebox(0,0)[lb]{\smash{\SetFigFont{12}{12.0}{it}carry}}}
\put(500,0){\makebox(0,0)[lb]{\smash{\SetFigFont{12}{12.0}{rm}[2,3]}}}
\put( 20,20){\line(0,1){ 40}}
\put( 20,20){\line(1,0){ 80}}
\put( 40,20){\framebox(60,40){}}
\put(150,20){\line(0,1){ 40}}
\put(150,20){\line(1,0){ 80}}
\put(210,20){\framebox(20,40){}}
\put( 20,0){\makebox(0,0)[lb]{\smash{\SetFigFont{12}{12.0}{it}OR}}}
\put( 55,0){\makebox(0,0)[lb]{\smash{\SetFigFont{12}{12.0}{rm}[1,2,3]}}}
\put(150,0){\makebox(0,0)[lb]{\smash{\SetFigFont{12}{12.0}{it}AND}}}
\put(200,0){\makebox(0,0)[lb]{\smash{\SetFigFont{12}{12.0}{rm}[3]}}}
\end{picture}

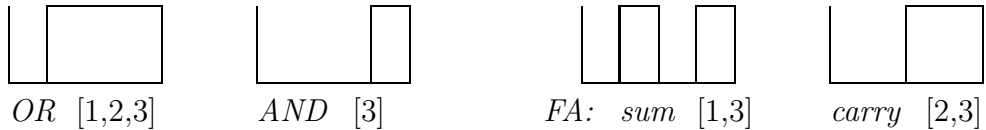
{\bf Fig 4}. ~OR, ~AND, ~Full Adder$(sum,carry)$

Symmetric functions $count$, typical for arithmetic. The well known
$OR$ function of $n$ inputs is symmetric, written $SF_n~[>$0]: at
least one high input, so only rank 0 is empty. The $n$-input $AND$
function is $SF_n~[n]$, active only if all $n$ inputs are high, so
only rank $n$ is full (containing just one minterm). And in a
3-input Full-Adder ($FA$): sum $s$=1 when 1 or 3 inputs are high,
so ranks [1,3] are full, written $s=SF[1,3]$, while the carry $c$=1
when 2 or 3 inputs are high, so $c=SF[2,3]$.

Most $BF$ however are not symmetric in all inputs, although many
have partial symmetries (in some inputs). A factored function $F$
cannot be symmetric, since inputs to different factors are not
equivalent. So an $SF$ has no factor, explaining why most logic
synthesis tools, based on factoring, have trouble with efficient
decomposition. \\This suggests putting $SF$'s in the Cell Library,
with $2^k~~SF_k$ cells of $k$ inputs, halving the number of cells
by using an inverter to exploit $SF(-X) = -SF(X)$.

{\bf $T$-cell library, threshold logic cells}

Threshold logic functions $TF<SF$ can implement any $SF$,
in a simple fashion.\\
{\it Def}: A {\bf threshold function} $T_k$ of $n$ inputs
has threshold $k \in [1,..,n]$ with $T_k$=1 whenever at
least $k$ inputs are active (high).

Any interval $[i, .. , j$-1] of $SF$ fullranks can be implemented
by the $AND$ of two threshold functions: ~$T_i~.~\overline{T_j}$.
So an $SF$ with $m$ fullrank intervals is the sum of $m ~TF$ pair
products.

For instance the FullAdder sum output (fig.4) with interval [1,2]
yields: $S[1,3]=(T_1.\overline{T_2})+T_3$, ~using the inverse of
carry $T_2$.

There are just $n$ $TF$ functions of $n$ inputs, with thresholds
$1,..,n$ - forming a compact and complete $T$-cell Library.
Including an inverter, a {\bf T-cell library} contains $sum(1,..,n)
= n(n+1)/2$ cells, that is 10 cells if $n$=4, or 15 cells for $n$=5.
~This is less than a complete {\bf S-cell library} of 1+(3+7+15)=26
cells ($n$=4), or 57 cells ($n$=5), which however will yield
smaller synthesized circuits (re section 6: further research).

\section{Planar cut and factoring}      

The two basic causes for {\bf asymmetry} are:~~~
  $factoring$ ~~and~~ $inverse$.

The smallest asymmetric functions are:~~~~~
   $a(b+c), ~~a+bc$ ~and~~ $\overline{a}~b, ~~\overline{a}+b$.

The first two cases use both (.) and (+) where the role of $a$
essentially differs from $b,c$ which are equivalent (permutable).
The last two cases are asymmetric in $(a,b)$, but symmetric
in $(\overline{a}, b)$. In general, {\bf input phasing} costs
little,  making a function more symmetric and increasing local
symmetries (with dense path sharing), essential for logic
optimization (fig.1,2)

{\bf Spectral product, and planar cut:~}
Function $F=G(X)~H(Y)$ is a $disjoint$ product if factors $G$ and
$H$ share no inputs, so $X \cap Y$ is empty. Multiplying the
rank spectra $sp(G)$ and $sp(H)$, as a convolution, yields the
spectrum of composition $F$:

~~~~~~~~~~$sp(F)=sp(G)*sp(H)$.

Order input sets $X$ and $Y$ adjacent in the gridplot. Then
this {\bf spectral product rule} follows since each path in
$G(X)$ is continued by (in product with) each path in $H(Y)$,
to form all paths (minterms) of length $|X|+|Y|$ in $F$. Let
$|X|=m$ then the gridplot of $F$ has diagonal $m$ consisting
of only planar nodes, with corresponding factor property:
$planar~cut$ (sect. 5 algorithm step 3).~~
Let $G=a~\#~b$ ~and~ $H=c+d+e$ with spectra $G[0,2,0]$ and
$H[0,3,3,1]$ then the product spectrum is $[0,3,3,1].[0,2,0]=
[0,0,6,6,2,0]$  by 'longhand' multiplication (without carry).

\section{'$Ortolog$' ~fast algorithm}          

The $Ortolog$ algorithm is designed for global yet fast detection
of (partial) symmetries, enhancing them by input phasing.
The rank spectrum is a simple and fast  symmetry test for
any sub function, by checking if each rank is full or empty.

The input format is that of a $PLA$ (2-level $or/and$ logic), hence
a list of $cubes$ as generalized minterms, each with all $n$ circuit
inputs (length $n$ strings over 1/0/- for input straight /inverse
/independent).
The algorithm is double recursive: start with a minimized 2-level
logic $BF_n(X)$ as a list of $m$ cubes, and proceed as follows:

\begin{enumerate}
\item $Core(a,b)$: for each input pair $(a,b)$ ~find the cubes
symmetric in $a,b$. \\Maximize each core by chosing input phase
$\overline{a}$ ~if~ $Core(\overline{a},b)$ has more cubes.

\item Input-expand maximal (phased) paircores to $Core(a,b,Y)$
with inputs $c$ (or $\overline{c}$) in rest input set $Y$.
Stop criterion: max $|Core| \times |inputs|^2$ prefers wide
(more inputs) over deep Core (more cubes). Select one
such 'best' multi input $Core(Z)$, symmetric for all
inputs in $Z \subseteq X$. Let $Y=\overline{Z}=X-Z$.

\item Factorize $Core(Z)$=$\sum_0^n G_r(Z)*H_r(Y)$ for ranks
  $r \leq n$ with non-zero symmetric rank- functions $G_r(Z)$
  as factors ($planar~cut$).

\item Recursively decompose (1-4) cofactors $H_r$ untill all
   components are symmetric.

\item Recursively decompose (1-5) remainder $F(X)-Core(Z)$,
  yielding an optimally phased network of symmetric functions
  coupled by inverters.
\end{enumerate}

{\bf Speedup option}: initially partition $F$ by collecting cubes
with equal number of dont-cares ($DC$ class), since cubes symmetric
in the same subset of inputs likely have the same number of $DC$'s.
Decompose the $k$ subfunc's $FDC_i$ separately:~ $F=\sum_1^k FDC_i$.

The $SF$ components can be implemented by $T$-cells, if a small
$T$-cell library is preferred. However, {\bf not decomposing} the $SF$
cells yields better area efficiency, using their grid plot as
layout pattern on silicon ($grid~template$), maximally sharing
logic paths.

The algorithm {\bf time complexity} is $O(n^2m)$, for a $BF_n$
list of $m$ cubes with $n$ inputs (step 1 is quadratic in $n$).
So only quadratic in the number of inputs (not exponential),
and linear in the number of cubes. This allows $very ~fast
~synthesis$ of many alternatives in a search for an optimal
binary code at a higher level: error correction codes in
Boolean circuit design [6][7][8] or state-machine logic:
$FSM$ state coding [9].

\subsection{Experiments}

The described symmetric synthesis with a cell
library of 15 $T$-cells (up to 5 inputs), was compared with a known
tool $Ambit$ (Cadence) using either a basic libary of $AND_n/OR_n/INV$
($n$=2..5) cells, or the usual extensive (full) libary of several
hundreds of cells. The logic {\it density} 'dens' is the filling \%
(non-$DC$) of the $PLA$ table to be decomposed. Rather than number
of cells, the total number of cell pitches (\#p) is compared in
Table 1, as area estimate:

\begin{verbatim}
     cct  inp cub dens  Synthesized  #pitches
    ------|--+--+  %   --Ambit--  Ortolog   Ratio
    binom5  6 32  74  (126)  128    148      0.86
    cordic 22 27  24  (135)  226    194      1.16
    table3 14 52  75  (448)  718    902      0.80
    parity  4  8 100  ( 18)   41     48      0.85
      Cell Library:  (Full)  AOI     TC     AOI/TC
\end{verbatim}
{\bf Table 1.} Synthesis areas (Standard cell \# pitches)

\section{Further research}     

{\bf Extend symmetric to planar functions}: The efficiency of
decomposing to a network of symmetric boolean functions clearly
depends on the amount of (local) symmetries in the initial $BF$.
Table 1 shows that restriction to a library of $AND/OR$ (column $AOI$)
resp. threshold $T$-cells (column $TC$) is too severe: results do not
compete with the usual large cell library, except the $cordic$ circuit
which has "much structure", viz. many local symmetries.

Symmetric components $SF_k$ (with dense sharing of logic paths)
should $not$ be mapped onto $T$-cells, but rather be implemented
directly as $planar$ compiled grid cells:

$Def$:~ a {\bf planar} Boolean function $PF_n$ has a planar grid-plot
(permute / invert inputs).

Notice that each symmetric $SF$ has only planar nodes in its
gridplot, hence is planar. Let a {\bf link} be a path of length=1
anywhere in a gridplot. Then any $SF_n$ is the 'template' for a
class of $PF_n$ easily derived from it by removing one or more
links. Obviously, any $PF_n$ has a unique smallest covering $SF_n$.

The class of $PF$ is much larger than $SF$, while being easily
derived by 'programming' (deleting links from) the $SF$'s as
templates. The number of links in any $SF_n$ is maximally
$\sum_1^n 2i = n(n+1)$, hence quadratic in $n$, rather than
exponential as in the case of look-up table $FPGA$'s.

The number of $PF_n$, between $|SF_n|=2^{n+1}$ and $|BF_n|=2^{2^n}$,
requires more research. All $BF_3$ are planar, and likely all
$BF_4$ as well, while non-planar $BF_n$ have $n \geq 5$.
\\[2ex]
{\large\bf Conclusions}

The symmetric T-cell library is too restricted to compete with the
usually very large cell libraries, since most $BF_n$ do not have many
sizable local symmetries. The area cost of lacking special cells
(e.g. $XOR$ in $parity$), and $T$-cell mapping of $SF's$ is high.

The $Ortolog$ algorithm performs $fast~global$ analysis, including
phase assignment, of local $BF_n$ symmetries. It detects and enhances,
by input phasing, the (dense) symmetric parts of a circuit, for
separate symmetric synthesis. The remaining (sparse) asymmetric logic
can be synthesized otherwise. ~~Flexible $compiled~cell$ logic synthesis,
using the larger class of {\bf planar} $BF$, can derive from symmetric
$SF_n$ as {\it programmable} $n \times n$ grid template.
\\[2ex]
{\large\bf References}

\begin{enumerate}
\item S. Akers: "Binary Decision Diagrams", IEEE Comp.
C-27, 509-516, June 1978

\item R. Bryant: "Graph-based algorithms for Boolean function
manipulation", IEEE Comp. C-35, 677-691, Aug 1986

\item L. Heinrich-Litan, P. Molitor: "Least Upper Bounds
for the size of OBDDs using Symmetry Principles",
IEEE Comp. C-49, 360-8, Apr 2000

\item J. v.Eijndhoven: "CMOS cell generation for Logic
Synthesis", proc. ASICON'94, 75-78, W.Y.Yuan (Ed)
Beijing, Oct 1994.

\item T. Courtney, et.al.: "Multiplexer based reconfiguration
for Virtex multipliers", {\it Field- Programmable Logic and
Applications}, FPL2000, 749-758, Villach, Austria, Aug 2000.

\item G. Muurling: "Fault tolerance in $IC$ design using error
correcting codes", \\MSc thesis TU-Delft, July 2000.

\item G. Muurling, et.al: "Error correction for combinational
logic circuits",  Benelux 21-st {\it Symposium on Information
Theory}, 25-31, Wassenaar, May 2000.

\item R. Kleihorst, N. Benschop: "Experiments with fault tolerant
$IC$ design using error correcting codes",  {\it International
Online Testing workshop}, Sicily, July 2001.

\item N. Benschop: "The structure of Constant Rank State Machines",
{\it Logic and Architecture Synthesis}, 167-176, G.Saucier (Ed.)
Paris, May 1990 (North-Holland, 1991)
\end{enumerate}

\end{document}